\documentclass[12pt]{article}
\usepackage{fullpage}
\usepackage{amsmath} 
\usepackage{amssymb}

\usepackage{theorem}
\newtheorem{theorem}{Theorem}[section]

\newtheorem{corollary}[theorem]{Corollary}

{\theorembodyfont{\rmfamily}

}
\def\lam{\lambda}
\def\alp{\alpha}
\def\kap{\kappa}
\def\sig{\sigma}
\def\gam{\gamma}
\def\os{\overline\sigma}
\def\l{\langle}
\def\r{\rangle}
\def\nek{,\ldots,}
\def\k{{\kappa_*^{+\gamma^{**}}}}

\title{On Some Configurations Related to
the Shelah Weak Hypothesis}
\author{Moti Gitik$^1$ and Saharon Shelah$^2$\\
\\
$^1$School of Mathematical Sciences\\
Tel Aviv University\\
Tel Aviv 69978, Israel\\
\\
$^2$Department of Mathematics\\
Hebrew University of Jerusalem\\
Jerusalem, Israel}
\date{}
\begin{document}
\thispagestyle{empty}
\maketitle
\begin{abstract}
We show that some cardinal arithmetic
configurations related to the negation of
the Shelah Weak Hypothesis and natural from
the forcing point of view are impossible.
\end{abstract}

\baselineskip=18pt
\section{Introduction}
The Shelah Weak Hypothesis (SWH) formulated
in [Sh:400A] states that for every cardinal
$\lam$ the number of singular cardinals
$\kappa < \lam$ with $pp\kappa \ge \lam$
is at most countable.  The negation of SWH
is one of the weakest statements on
cardinal arithmetic whose consistency is
unknown.  Clearly, SWH follows from GCH or
even from the Shelah Strong Hypothesis
saying that for every singular $\kappa$,
$pp\kappa = \kappa^+$.  On the other hand
by [Sh-g], ``$|pcfa| > |a|$" implies
$\neg$SWH.
The forcing construction of [Gi-Sh] and
[Gi-Ma] show that any finite or countable
number of $\kappa$'s with $pp\kappa >
\lam$ is consistently possible.

The present paper grew from an attempt made
by the first author to force $\neg$SWH
using a forcing of type of [Gi].  One of the
features of this forcing is that it does
not add new bounded subsets to a cardinal
while blowing its power.  Here we show (in
ZFC) that some configurations which are
very natural from the forcing point of view
are just impossible.

The first theorem under stronger
assumptions was proved by the first author;
the second author was able to weaken
the assumptions and find  a more elegant
proof.  Most of the generalizations are due
to the second author.  The second theorem
is due solely to the second author.
\section{Main Results}

\medskip

\noindent
{\bf Theorem 1}
{\it The following is impossible:
\begin{itemize}
\item[{\rm (a)}] $\kappa_1 < \kappa_*$, $cf \kappa_1 =
\aleph_0$, $cf\kappa_* > 2^{\aleph_0}$
\item[{\rm (b)}] for every large enough $\mu <
\kappa_1$ of cofinality $(2^{\aleph_0})^+$
we have $pp\mu = \mu^+$
\item[{\rm (c)}] $\kappa_* = \sup \{ \mu \mid \mu <
\kappa_*, cf \mu = \aleph_0\ {\rm and}\
pp\mu > \kappa_*^+\}$
\item[{\rm (d)}] there are a strictly increasing
sequence $\langle \lambda_\alpha \mid
\alpha < cf \kappa_*\rangle$ of regular
cardinals between $\kappa_1$ and $\kappa_*$
unbounded in $\kappa_*$, a filter $D$ on
$\omega$ containing all cofinite subsets of
$\omega$ and a sequence of functions
$\langle f_{\lambda_\alpha} \mid \alpha <
cf \kappa_*\rangle$ such that
\begin{itemize}
\item[$(\alpha)$] $f_{\lambda_\alpha}:\omega
\to Reg \cap \kappa_1\backslash
(2^{\aleph_0})^+$
\item[$(\beta)$] $\lim_D f_{\lambda_\alpha} =
\kappa_1$
\item[$(\gamma)$] $\lambda_\alpha = tcf
(\prod_{n < \omega} f_{\lambda_\alpha} (n)
/D)$
\item[$(\delta)$] $\alpha < \beta < cf
\kappa_*$ implies $f_{\lambda_\alpha} <_D
f_{\lambda_\beta}$
\item[$(\epsilon)$] if $\alpha < \beta < cf
\kappa_*$ and $\lambda \in Reg \cap
\lambda_\beta \backslash \lambda_\alpha^+$
then there is a function $f_\lambda:\omega
\to Reg \cap \kappa_1 \backslash
(2^{\aleph_0})^+$ such that
$f_{\lambda_\alpha} <_D f_\lambda <_D
f_{\lambda_\beta}$ and $\lambda = tcf
(\prod_{n<\omega} f_\lambda(n)/D)$.
\end{itemize}
\end{itemize}
  }
 
\medskip
\noindent
{\bf Discussion 1.1} 
\begin{itemize}
\item[(1)] The assumption (c)
is a form of $\neg$SWH which claims that
there are more than $2^{\aleph_0}$ singular
cardinals of cofinality $\aleph_0$ with
$pp$ above their supremum.
\item[(2)] The assumption (d) holds
naturally in forcing constructions with $D
=$ the filter of cofinite subsets of
$\omega$. But it
seems to be problematic in ZFC.  In [Sh-g,
II\S1] proof of a weak relative is a major
result.
\item[(3)] See [Sh-g, VI] for uncountable
cofinalities.
\end{itemize}
\medskip

\noindent
{\bf Proof.} Suppose otherwise.  W. l. of
g. we can assume that $cf \kappa_* =
(2^{\aleph_0})^+$.  Also, replacing $\l
\lambda_\alpha \mid \alpha <
(2^{\aleph_0})^+ \rangle$ by its
restriction to an unbounded subset, we can
assume that the following holds (see [Sh-g]
for this type of argument)\\[.5em]
$(*)$ for every $n < \omega$,
$\langle f_{\lambda_\alpha} (n)\mid \alpha
< (2^{\aleph_0})^+\rangle$ is strictly
increasing and, if $f_* (n) =
\bigcup_{\alpha < (2^{\aleph_0})^+}
f_{\lambda_\alpha} (n)$ then $f_* (n) <
f_{\lambda_0} (n+1)$.

Now for every $\alpha < (2^{\aleph_0})^+$
and $\lambda \in Reg \cap
\lambda_{\alpha+1} \backslash
\lambda_\alpha$ we use $(\epsilon)$ and
find a function $f_\lambda:\omega \to Reg
\cap \kappa_1\backslash (2^{\aleph_0})^+$
such that $\lambda = tcf (\prod_{n<\omega}
f_\lambda(n) / D)$ and for every
$n <\omega$, $f_{\lambda_\alpha} (n) <
f_\lambda(n) < f_{\lambda_{\alpha+1}} (n)$.

Clearly, $\l f_* (n) \mid n <
\omega\rangle$ is strictly increasing with
limit $\kappa_1$ and $cf (f_* (n)) =
(2^{\aleph_0})^+$ for every $n < \omega$.
Using (b), we can assume removing finitely
many $n$'s, if necessary,  that $pp(f_*(n)) =
(f_*(n))^+$ for every $n < \omega$.
Let $D_*$ by an ultrafilter on $\omega$
extending $D$.
Let $\mu_* = tcf \prod_{n<\omega} ((f_*
(n))^+/D_*)$.  It is well defined since
$D_*$ is an ultrafilter.  By (c), w.l. of
g., for every $\alpha < (2^{\aleph_0})^+$
there is $\kappa_\alp$, $\lam_\alp <
\kap_\alp < \lam_{\alp+1}$, $cf \kap_\alp =
\aleph_0$ and $pp\kap_\alp \ge
\kap_*^{++}$. Hence, there are
$\tau^2_{\alp,n} \in Reg \cap \kap_\alp
\backslash \lam_\alp^{++}$ $(n < \omega)$
and a filter $D_\alp$ on $\omega$
continuing all cofinite subsets of $\omega$
such that $\kap_*^{++} = tcf (\prod_{n <
\omega} \tau^2_{\alp,n} /
D_\alp)$.  By [Sh-g], we can then find
$\tau^1_{\alp,n} \in Reg \cap
\tau^2_{\alp,n}\backslash \lam_\alp^+$ such
that $\kap_*^+ = tcf (\prod_{n<\omega}
\tau^1_{\alp,n}/D_\alp)$ (note
that we are doing this separately for each
$\alp < (2^{\aleph_0})^+)$.  Let
$a_\alp^{m,\ell}= \{f_{\tau^\ell_{\alp,n}}
(m)\mid n < \omega\}$ for every $m <
\omega$ and $\ell \in \{1,2\}$.  Set 
$a^m_\alp = a_\alp^{m,1} \cup
a_\alp^{m,2}$, $a^m = \bigcup_{\alp <
(2^{\aleph_0})^+} a_\alp^m$ and $a =
\bigcup_{m < \omega} a^m$.  All these sets
consists of regular cardinals above
$(2^{\aleph_0})^+$, $a_\alp^m$'s are
countable, $a^m$'s and $a$ have cardinality
of at most $(2^{\aleph_0})^+$.  Also
$a_\alpha^m \subseteq [f_{\lambda_\alpha}
(m), f_{\lambda_\alp} (m + 1))$.  Clearly,
$a^m$ $(m < \omega)$ is an unbounded subset
of $f_* (m) \cap Reg$ of order type
$(2^{\aleph_0})^+$, since $\langle
f_{\lam_\alp} (m) \mid \alp <
(2^{\aleph_0})^+\rangle$ is increasing with
limit $f_*(m)$. Then, $(f_* (m))^+ = tcf
(\prod a^m / J_{a^m}^{{\rm
bounded}})$, as $pp (f_*(m)) = (f_*
(m))^+$.  Hence, by [Sh-g], $(f_*(m))^+ \in
pcf a^m \subseteq  pcf a$, for every $m <
\omega$. Again, by [Sh-g], $pcf (\{(f_*
(m))^+ \mid m < \omega\}) \subseteq pcfa$.
But $\mu_* = tcf (\prod_{n < \omega} (f_*
(n))^+/ D_*)$, hence $\mu_* \in
pcf a$.  Let $\langle b_\sigma \mid \sigma
\in pcf a\rangle$ be a generating sequence
for $a$ (see [Sh-g] or [Sh:506]).  W.l. of g.
if $\mu_*\ne \kappa_*^{+\ell}$ for $\ell
\in \{1,2\}$, then $b _{\mu_*} \cap
b_{\kappa_*^{+\ell}} = \emptyset$.  Let
$\ell^* \in \{1,2\}$ be such that $\mu_*
\ne \kap_*^{+\ell^*}$.
\medskip

\noindent
{\bf Claim 1.2} {\it The set $A = \{m <
\omega\mid \ \hbox{for some}\ \alp <
(2^{\aleph_0})^+$, $\bigcup_{\beta\in
[\alpha,(2^{\aleph_0})^+)} a^m_\beta
\subseteq b_{\mu_*}\}$ is in $D_*$.
}
\medskip

\noindent
{\bf Proof.} Otherwise $\omega \backslash A
\in D_*$ and for $m \in \omega \backslash
A$, $f_* (m) = \sup (a^m \backslash
b_{\mu_*})$.  Hence $(f_* (m))^+ \in pcf (a
\backslash b_{\mu_*})$.  So, $pcf (\{(f_*
(m))^+ \mid m\in \omega \backslash A\})
\subseteq pcf (a \backslash b_{\mu_*})$.
But $\omega\backslash A \in D_*$ and $\mu_*
= tcf (\prod_{m < \omega} (f_*
(m))^+/ D_*)$.  Hence $\mu_* \in
pcf (a \backslash b_{\mu_*})$.
Contradicting the choice of $b_{\mu_*}$.\\
$\square$ of the claim.

For $m \in A$ let $\alp_m$ be the minimal
$\alp$ such that $\bigcup_{\beta \in
[\alp,(2^{\alp_0})^+)} a^m_\beta \subseteq
b_{\mu_*}$.  Set $\alp_* = \bigcup_{m\in A}
\alp_m$.  Clearly, $\alp_* <
(2^{\aleph_0})^+$.  Let $a' = \bigcup \{
a^m_\beta \mid m \in A, \beta \in [\alp_*,
(2^{\aleph_0})^+)\}$. Then $a' \subseteq
b_{\mu_*}$ and hence $\kappa_*^{+\ell^*}
\not\in pcf a'$.
However, $m \in A$ and $n < \omega$ imply
that $f_{\tau^{\ell^*}_{\alp_*,n}} (m) \in
a_{\alp_*}^{m,\ell^*} \subseteq
a_{\alp_*}^m \subseteq a'$.
So, for each $n < \omega$ we have
$$\{f_{\tau^{\ell^*}_{\alp_*,n}} (m) \mid m
\in A\} \subseteq a'\ .$$
Hence $pcf \{f_{\tau^{\ell^*}_{\alp_*, n}}
(m) \mid m \in A\} \subseteq pcf a'$.  But
as $A \in D_*$, $\tau^{\ell^*}_{\alp_*,n} =
tcf (\prod_{m\in A}
f_{\tau^{\ell^*}_{\alp_*,n}} (m)/
D)$. So, for every $n < \omega$,
$\tau_{\alp_*,n}^{\ell^*} \in pcf a'$.  Then
by [Sh-g],  $pcf \{
\tau^{\ell^*}_{\alp^*,n}\mid n < \omega\}
\subseteq pcf a'$.  But $\kappa_*^{+\ell^*}
= tcf (\prod_{n < \omega}
\tau_{\alp^*,n}^{\ell^*} / D_{\alp^*})$. So,
$\kappa_*^{+\ell^*} \in pcf a'$.
Contradiction.\\
$\square$

\noindent
{\bf Remark 1.3}
\begin{itemize}
\item[(1)] We can replace ``$cf \kappa_*>
2^{\aleph_0}$" by ``$cf \kappa_* >
\aleph_0$" provided that (d) of the theorem
is strengthened by adding the condition $(*)$
introduced in the beginning of the proof
and $(2^\aleph)^+$ is replaced by $\aleph_1$
in (b). 
\item[(2)] It is possible to weaken
``$pp\mu > \kappa^+_*$" in (c) of the
theorem to ``$pp\mu \ge \kappa_*$",
replacing $(2^{\aleph_0})^+$ in (b) by
$\aleph_1$.  Just after $(*)$ is obtained
using $cf \kappa_* \ge (2^{\aleph_0})^+$,
we can replace $\kappa_*, \kappa_*^+,
\kappa_*^{++}$ by the limit of first
$\aleph_1$, $\lambda_\alp$'s, its successor
and its double successor.  This is provided
that for every $\alp < \omega_1$ there is
$\kappa_\alp,$
$\lam_\alp<\kappa_\alp<\lam_{\alp + 1}$
with $pp \kappa_\alp \ge \lam_*^{++}$,
where $\lam_* = \bigcup_{\alp<\omega_1}
\lam_\alp$.  The condition ``$pp\mu \ge
\kappa_*$" can be used easily to construct
such $\langle \lam_\alpha\mid \alp <
\omega_1\rangle$.
\item[(3)] It is possible to replace ``$cf
\kappa_* > 2^{\aleph_0}$" by ``$\forall
\alpha < cf \kappa_*$ $(|\alp |^{\aleph_0}
< \kappa_*)$" just use $cf \kappa_*$
instead of $(2^{\aleph_0})^+$ in the
proof.
\end{itemize}

The following is parallel to Solovay's
result that SCH holds above a strongly
compact cardinal.
\setcounter{section}{1}
\setcounter{theorem}{3}
\begin{corollary}
Suppose that the following holds: $\kappa$
is a cardinal such that 
\begin{itemize}
\item[{\rm (a)}] for any given
cardinal $\lam$ it is possible to force
$2^\kappa \ge \lam$ by $\kappa^{++}$-c.c.
forcing not adding new bounded subsets to
$\kappa$ and adding $\lam$
$\omega$-sequences $\langle f_\alp \mid
\alp < \lam\rangle$ to $\kappa$ such that
$\alp < \beta \to f_\alp < f_\beta$ (mod
finite) and $\delta \in (\kappa, \lam]$
regular cardinal implies that $f_\delta(n)$
is regular cardinal for every $n < \omega$.
\item[{\rm (b)}] $pp(\mu) = \mu^+$ for
every large enough $\mu < \kappa$ of
cofinality $\aleph_1$.
\\Then above $\kappa$ the following version
of SWH holds:\\ for every cardinal $\lam$ the
set $\{\mu \mid \kappa < \mu < \lam, cf \mu
= \aleph_0, pp\mu > \lam^+\}$ is at most
countable.
\end{itemize}
\end{corollary}
\medskip

\noindent
{\bf Remark.} Forcing notion of [Gi-Ma] or
[Gi] satisfy (a).
\medskip

\noindent
{\bf Proof.} Suppose otherwise.  Let
$\kappa_*$ be the first cardinal such that
the set $\{\mu \mid \kappa < \mu <
\kappa_*, cf \mu = \aleph_0,pp\mu >
\kappa_*^+\}$ is uncountable.  Clearly,
$cf\kappa_* =\aleph_1$.  Now we force with
the forcing of (a) and make $2^\kappa \ge
\kappa_*$. The $\omega$-sequences produced
by such forcing will satisfy $(*)$ of the
proof of Theorem 1 with $D$ equal to the
filter of cofinite sets.  The chain
condition of the forcing insures that the
cardinal arithmetic does not change above
$\kappa$.  No new bounded subsets are added
to $\kappa$, hence  (b) of the
statement of the corollary still holds.  Now
Theorem 1 (actually using 1.3(2)) provides
a contradiction.\\
$\square$	
\medskip

Repeating the proof of Theorem 1 we can
show the following generalization:
\begin{theorem}
The following is impossible:
\begin{itemize}
\item[{\rm (a)}] $\kappa_1<\kappa_*$, $cf
\kappa_1 = \aleph_0$, $cf \kappa_* >
2^{\aleph_0}$
\item[{\rm (b)}] there is $\ell$, $1 \le \ell <
\omega$ such that for every $\mu < \kappa_1$
of cofinality $(2^{\aleph_0})^+$ we have
$pp\mu \le \mu^{+\ell}$
\item[{\rm (c)}] $\kappa_* = \sup \{ \mu \mid \mu
< \kappa_*, cf \mu = \aleph_0\ {\rm and}\
pp\mu > \kappa_*^{+\ell}\}$
\item[{\rm (d)}] the same as in Theorem 1.
\end{itemize}
\end{theorem}
If we allow infinite gaps between $\mu$ and
$pp\mu$ in (b) of 1.5, then the following
can be shown:
\begin{theorem}
Assume that 
\begin{itemize}
\item[{\rm (a)}] $\kappa_1 < \kappa_*$,
$cf \kappa_1 = \aleph_0$, $cf \kappa_* =
\theta > \aleph_0$, $\alp^* < \kappa_1$,
$cf \alpha^* > \aleph_0$
\item[{\rm (b)}] for every large enough
$\mu < \kappa_1$ of cofinality $\theta$ we
have $pp(\mu) < \mu^{+ \alp^*}$ 
\item[{\rm (c)}] for some $\beta^*$,
$\kappa_* = \sup \{ \mu\mid \mu < \kappa_*,
cf \mu = \aleph_0 \ {\rm and}\  pp\mu \ge
\kappa_*^{+\beta^*}\}$
\item[{\rm (d)}] the condition (d) of
Theorem 1 and $(*)$ of its proof
\end{itemize}
Then $\beta^* < \sigma^{+4}$ for some
$\sigma < \alpha^*$.
\end{theorem}

\medskip

\noindent
{\bf Sketch of the proof.}
Suppose otherwise.  We define $f_* (n)$'s
as in Theorem 1.  Now $cf f_* (n) = \theta$
and so $pp (f_* (n)) < (f_* (n))^{+\alp^*}$
for every $n < \omega$. Find $\sigma <
\alp^*$ such that for every $n < \omega$,
 $pp(f^*(n)) \le (f_*
(n))^{+\sigma}$.  Here we use that $cf
\alp^* > \aleph_0$.  Instead of one
$\mu_*$ in the proof of Theorem 1 (or
finitely many cardinals in 1.6) we consider
$pcf \{(f_* (n))^{+\sigma'} \mid n <
\omega,\sigma'\le \sigma\} \cap (\kappa_*,
\kappa_*^{+\beta^*}]$. By the assumption we
made, $\beta^* \ge \sigma^{+4}$. Then there
should be $\kappa_*^{+\ell^*} \not\in pcf
\{(f_* (n))^{+\sigma'} \mid n < \omega,
\sigma' \le \sigma\}$ for some $\ell^*$, $1
\le \ell^* \le \beta^*$.  This follows by
results of [Sh:g, IX], see also [Sh:g, E12,
4.18 (b)]. The rest of the proof is as
those of Theorem 1, only we use [Sh:g, I,
3.2(5)] to include $pcf\{(f_*
(n))^{+\sigma'}\mid n<\omega, \sigma' \le
\sigma\}$ into a union of finitely many
$pcf$-generators.\\
$\square$

Now we turn to another theorem which
provides a different proof of Theorem 1 and
some of its generalizations.

\medskip

\noindent
{\bf Theorem 2} {\it Suppose that
\begin{itemize}
\item[{\rm (a)}] $\kappa_0 < \kappa_1 <
\kappa_*$, $1\le n^*< \omega$, $n^* <
\gamma^* < \theta$ and $\gamma^*$ is a
successor ordinal
\item[{\rm (b)}] $\theta = cfd \theta <
\kappa_0$ and for every $\alpha < \theta$,
$|\alpha|^{\aleph_0} < \theta$
\item[{\rm (c)}] $cf \kappa_1 = \aleph_0$
and $pp(\kappa_1) \ge \kappa_*^{+\gamma^*}$
\item[{\rm (d)}] if $\mu \in
(\kappa_0,\kappa_1)$ and $cf \mu = \theta$
then $pp\mu \le \mu^{+n^*}$.
\end{itemize}
\underline{Then} the following holds 
\begin{itemize}
\item[{\rm (1)}] For every nonprincipal ultrafilter $D$
on $\omega$ and a sequence $\overline
\sigma^* = \langle \sigma^*_\ell \mid \ell
< \omega \rangle$ with $\kappa_1 = \lim_D
\overline \sigma^*$ and $\sigma^*_\ell$
$(\ell < \omega)$ a limit cardinal of
cofinality $\ge \theta$ in the interval
$(\kappa_0,\kappa_1)$ \underline{there are} a set
$w \subseteq \gamma^* + 1$ consisting of at
most $n^*$ elements and a sequence
$\overline \sigma^{**} = \langle
\sigma_\ell^{**} \mid \ell <
\omega\rangle$, $\kappa_0 <
\sigma_\ell^{**}< \sigma_\ell^*$\enskip $(\ell <
\omega)$ \underline{such that}
\\[.5em]
$(*)_1$  if $a \in [R_{D,\overline
\sigma^*, \overline
\sigma^{**}}]^{\aleph_0}$, $\beta \le
\gamma^*$ and $\kappa_*^{+\beta}\in pcfa$
then $\beta \in w$, \\
where
$R_{D,\overline \sigma^*, \overline
\sigma^{**}} = \{ tcf (\Pi \overline \sigma
/D)\mid \overline \sigma = \langle \sigma_n
\mid n < \omega \rangle,\ \sigma^{**}_n
\le \sigma_n = cf \sigma_n < \sigma_n^*
\enskip (n<\omega)\} \cap [\kappa_1, \kappa_*)$.
\item[{\rm (2)}] There are $\alpha^* < \theta$ and a
sequence $\langle R_\alpha \mid \alpha <
\alpha^*\rangle$ with $\bigcup_{\alpha <
\alpha^*} R_\alpha = Reg \cap \kappa_*
\backslash \kappa_1$ \underline{so that}
\\[.5em]
$(*)_2$ for every $\alpha < \alpha^*$ there
is $w \le \gamma^* + 1$ consisting of at
most $n^*$ elements such that\\
if $a \in [R_\alpha]^{\aleph_0}$, $\beta
\le \gamma^*$ and $\kappa_*^{+\beta} \in
pcf a$ then $\beta \in w$.
\item[{\rm (3)}] Let $D$ be a nonprincipal ultrafilter
on $\omega$. There is a partition $\langle
I_\rho \mid \rho < \rho^*\rangle$, $\rho^*
< \theta$ of $Reg \cap \kappa_1 \backslash
\kappa_0$ into closed open intervals (i.e.
of the form $[x,y)$) with $\langle \min
I_\rho \mid \rho < \rho^*\rangle$ strictly
increasing such that
\\[.5em]
$(*)_3$ for every sequence $\langle \rho_n
\mid n < \omega \rangle$ of ordinals below
$\rho^*$ with $\lim_D \langle\min
I_{\rho_n} \mid n < \omega \rangle =
\kappa_1$
$$\big\{tcf (\prod_{n < \omega}\sigma_n
/D)\mid \sigma_n \in I_{\rho_n}\ {\rm
for}\ n < \omega\}\cap [\kappa_1,
\kappa_*)$$
is included in one of $R_\alpha$'s ($\alpha
< \alpha^*)$ from a sequence $\langle
R_\alpha \mid \alpha < \alpha^*\rangle$
$(\alpha^* < \theta)$ satisfying $(*)_2$.
\end{itemize}

}
\medskip

\noindent
{\bf Remark 2.1} Part (1) is close to [Sh:g, IX 1.x].
\medskip

\noindent
{\bf Proof of (2) and (3) from (1)}. As
$\kappa_*^{+\gamma^*} \le pp(\kappa_1)$
there are a countable unbounded $a
\subseteq \kappa_1 \cap Reg \backslash
\kappa_0$ and an ultrafilter $D_0$ on $a$
containing all cobounded subsets of $a$ with
$\kappa_*^{+\gamma^*} = tcf (\Pi a/ D_0)$.
Let $a= \{ \lambda_n \mid n < \omega\}$ and
$D = \{ A \subseteq \omega \mid \{
\lambda_n \mid n \in A\} \in D_0\}$.  Now,
by [Sh:g, II], for every regular $\tau \in
\kappa_*^{+\gamma^*} \backslash \kappa_1$
we can find $\overline \sigma = \langle
\sigma_n \mid n < \omega \rangle$,
$\sigma_n \in Reg \cap \kappa_1 \backslash
\kappa_0$ $(n < \omega)$, $\lim_D \overline
\sigma = \kappa_1$ such that $\tau = tcf
(\Pi \overline \sigma / D)$.

Fix $\chi$ to be a large enough cardinal.  Let
$M \prec (H(\chi), \epsilon)$ be such that $|M|<
\theta$, $^\omega M \subseteq M$,
$\{\kappa_0,\kappa_1, \theta, D, \kappa_*\}
\in M$ and $M \cap \theta \in \theta$.
There is such $M$ since we assumed (b).
Consider the following set $\Phi = \{
\overline \sigma^* \mid \overline \sigma^*
= \langle \sigma^*_n \mid n < \omega
\rangle$, $\lim_D \overline \sigma^* =
\kappa_1$ and for every $n < \omega$,
$\sigma_n^* \in M \cap [\kappa_0^+,
\kappa_1)$ is a limit cardinal of
cofinality $\ge \theta\}$. Clearly, $\Phi
\subseteq M$ since ${}^\omega\! M \subseteq M$.
Now, by (1), applied with $D$ defined above
for each $\overline \sigma^* \in \Phi$
there will be $\overline \sigma^{**}$ for
which $(*)_1$ holds.  By elementarity,
there is such $\overline \sigma^{**}$ in
$M$. Denote it by $\overline \sigma^{**}
[\overline \sigma^*]$.  Define $\langle
R_\alpha \mid \alpha < \alpha^*\rangle$ to
be an enumeration of the set
$\{R_{D,\overline
\sigma^*,\overline\sigma^{**}[\overline
\sigma^*]} \mid \overline \sigma^* \in
\Phi\} \cup \{\{ tcf (\prod_{n < \omega}
\sigma_n / D)\}\mid \sigma_n \in M \cap
\kappa_1 \cap Reg \backslash \kappa_0 \
{\rm and}\ \lim\limits_{n < \omega \ D } \sigma_n
= \kappa_1 \}$. Then $\alpha^* < \theta$
since $|M| < \theta$.  Clearly here $(*)_1$
implies $(*)_2$.  So, in order to complete
the proof of (2) it remains to show that
$Reg \cap \kappa_* \backslash \kappa_1 =
\bigcup_{\alpha < \alpha^*} R_\alpha$.
Let $\tau \in Reg \cap \kappa_* \backslash
\kappa_1$.  Then for some $\overline \sigma
= \langle \sigma_n\mid n < \omega \rangle$,
$\sigma_n \in Reg \cap \kappa_1 \backslash
\kappa_0$ $(n < \omega)$, $\lim_D \overline
\sigma = \kappa_1$, $\tau = tcf (\Pi
\overline \sigma / D)$.  Let $A = \{ n <
\omega \mid \sigma_n \in M\}$.
\medskip

\noindent
{\bf Case 1.} $A \in D$.

 Then, w.l. of g.
we can assume that $A = \omega$ (just if
$\sigma_n \not\in M$ replace it by
$\kappa_0^+$).  But then $\tau$ appears in
the second part of the union defining
$\langle R_\alpha\mid \alpha<
\alpha^*\rangle$.
\medskip

\noindent
{\bf Case 2.} $A \not \in D$.
 
Clearly $\kappa_1 \ge \kappa_0^{+\theta}$,
since otherwise $\kappa_1 \cap Reg
\subseteq M$ and Case 2 cannot occur.  So
w.l. of g. we can assume that $A =
\emptyset$. Let for $n < \omega$,
$\sigma_n^* = \min (M \cap \kappa_1
\backslash \sigma_n)$.  Such $\sigma_n^*$
is well defined since $\kappa_1 \in M$,
$cf\kappa_1 = \aleph_0$ and hence $\kappa_1
= \sup (\kappa_1 \cap M)$.  Also,
$\sigma_n^*$ should be a limit cardinal of
cofinality $\ge \theta$ as $M \cap \theta
\in \theta$.  So $\overline \sigma^* =
\langle \sigma_n^* \mid n < \omega \rangle
\in \Phi$. Let $\overline \sigma^{**} =
\overline \sigma^{**}[\overline \sigma^*]$.
Now, for every $n < \omega$, $\kappa_0^+
\le \sigma_n^{**} < \sigma_n^*$ and
$\sigma_n^{**} \in M$.  Hence,
$\sigma_n^{**} < \sigma_n < \sigma_n^*$ for
every $n < \omega$.  Then $tcf (\Pi
\overline \sigma/D) = \tau \in
R_{D,\overline \sigma^*,\overline \sigma^{**}}$
by $(*)_1$ and we are done.

This completes the proof of (2) from (1).

Let us turn now to (3).  Here we are given
a nonprincipal ultrafilter $D$.  Define $M$
and $\langle R_\alp \mid \alp <
\alp^*\rangle$ as above using this $D$.
For every $\nu \in M \cap \kap_1 \backslash
\kap_0$ a limit cardinal of cofinality $\ge
\theta$ denote $\sup (M \cap \nu)$ by
$\nu(M)$.  Let $\langle I_\rho \mid \rho <
\rho^*\rangle$ be the increasing
enumeration of the following disjoint
intervals:\\
$\{Reg \cap [\nu(M), \nu] \mid \nu \in M
\cap \kap_1$ is a limit cardinal of
cofinality $\ge \theta\}\cup \{ \{
\nu\}\mid \nu \in M,cf \nu = \nu \}$.

Clearly, $\rho^* < \theta$, since $|M| <
\theta$. Let us check that $(*)_3$ holds. 
So let $\langle \rho_n \mid n <
\omega\rangle$ be a sequence of ordinals
below $\rho^*$ with $\lim_D \langle \min
I_{\rho_n} \mid n < \omega\rangle =
\kappa_1$ and let $\sigma_n \in I_{\rho_n}$
for $n < \omega$.  Consider $\tau = tcf
(\prod_{n<\omega} \sigma_n / D)$.  Let $A =
\{ n < \omega \mid \sigma_n \in M\}$.  As
above we can concentrate on the situation
when $A = \emptyset$ (i.e. Case 2).  Define
$\overline \sigma^*$ and $\overline
\sigma^{**}$ as in Case 2.
Then for every $n < \omega$, $\sigma_n^{**}
< \sigma_n^*$ and $\sigma_n^{**} \in M$.
But $\sigma_n^* = \min (M \cap \kappa_1
\backslash \sigma_n)$ is a limit cardinal
of cofinality $\ge \theta$ in $M$. Let
$\tilde \rho_n$ denote the left side of
the interval $I_{\rho_n}$.  Then
$\sigma_n^* = \tilde \rho_n$, since $\tilde
\rho_n \in M$ is a limit cardinal of
cofinality $\ge \theta$ and $\sigma_n \in
I_{\rho_n} = (\sup (M \cap \tilde \rho_n),
\tilde \rho_n) \cap Reg$.  Also the last
equality implies that $\sigma_n >
\sigma_n^{**}$.  Then $\tau = tcf (\prod_{n
< \omega} \sigma_n /D) \in R_{D, \overline
\sigma^*, \overline \sigma^{**}}$ and we
are done.
\medskip

\noindent{\bf Proof of (1).}
Suppose otherwise.  Let $D$ be a
nonprincipal ultrafilter on $\omega$ and 
$\os^* = \l \sig^*_n \mid n < \omega \r$ a
sequence of limit cardinals of cofinality
$\ge \theta$ in the interval
$(\kap_0,\kap_1)$ with $\kap_1 = \lim_D
\os^*$ witnessing the failure of (1).  We
choose by induction on $\xi < \theta$
cardinals $\sig_{\xi,n},
\tau^k_\xi,\sig^k_{\xi,n}$ $(n,k<\omega)$
so that
\begin{itemize}
\item[($\alp$)] $\kap_0^+ \le \sig_{\xi,n}
< \sig^*_n$
\item[($\beta$)] $\xi < \xi'$ implies
$\sig_{\xi,n}^k < \sig_{\xi',n}$
\item[($\gamma)$] $\tau^k_\xi \in Reg \cap
\kappa_* \backslash \kappa_1$
\item[$(\delta)$] $\kap_*^{+\gamma} \cap
pcf (\{ \tau^k_\xi \mid k < \omega\})
\backslash \kap_*$
has at least $n^* + 1$ members
\item[$(\epsilon)$] $\sig_{\xi,n} <
\sig^k_{\xi,n} < \sig_n^*$ and
$\sig_{\xi,n}^k$ is regular
\item[$(\xi)$] $tcf (\prod_{n<\omega}
\sig^k_{\xi,n} / D) = \tau_\xi^k$
\item[$(\eta)$] $\xi < \xi'$ implies that
$\sig_{\xi,n} < \sig_{\xi',n}$.
\end{itemize}

In order to carry out the construction we
choose first at stage $\xi$, $\sig_{\xi,n}$
satisfying $(\alp),(\beta)$.  This is
possible, since $\sigma_n^*$ is a limit
cardinal $> \kap_0$ of cofinality $\ge
\theta$.  Second, as $\l \sig_{\xi,n} \mid
n < \omega\r$ cannot serve as $\os^{**}$ in
$(*)_1$ by our assumption, there are
$\tau_\xi^k \in R_{D,\os^*,\l
\sig_{\xi,n}\mid n < \omega\r}$ for $k <
\omega$ such that $pcf (\{ \tau^k_\xi \mid
k < \omega\}) \cap (\kap_*, \kap_*^{+\gamma^*}]$
has at least $n^* + 1$ members.  So clauses
$(\gamma),(\delta)$ hold. By the definition
of $R_{D,\os^*,\l\sig_{\xi,n} \mid n <
\omega\r}$, we can find for each $k <
\omega$, 
$\sigma_{\xi,n}^k \in Reg \cap \sigma_n^*
\backslash \sigma_{\xi, n}$ such that $tcf
(\prod_{n < \omega} \sigma^k_{\xi,n}, D) =
\tau_\xi^k$.  So clauses $(\epsilon)$
and $(\xi)$ hold.  The clause $(\eta)$ is
implied by the previous ones.  So, we have
finished the inductive construction.

Now, for every $n < \omega$, as $\l
\sig_{\xi,n} \mid \xi < \theta\r$ is
strictly increasing, its limit $\sig_n =
\bigcup_{\xi<\theta} \sig_{\xi,n}$ is a
singular cardinal of cofinality $\theta$.
Also, clearly, $\sig_n \in [\kap_0^+,\kap_1)$.
Hence, by the assumption $(d)$ of the
theorem, $pp\sig_n \le \sig^{+n^*}_n$. For
$\ell = 1 \nek n^*$ let $\lam_\ell = tcf
(\prod_{n < \omega} \sig_n^{+\ell}/D)$.
Set $w^* = \{ \alpha \le \gamma^* \mid
\kappa_*^{+\alpha} = \lam_\ell$ for some
$\ell$, $1 \le \ell \le n^*\}$.  Then
$w^*$ is a set of $\le n^*$ ordinals
below $\gamma^* + 1$. Let $a_n = \{
\sigma_{\xi,n}^k \mid k < \omega, \xi <
\theta\}$ and $a = \bigcup_{n < \omega} a_n
\cup \{ \sig_n^{+\ell} \mid n < \omega, 1
\le \ell \le n^*\}$.  So, $a$ is a set of
$\le \theta < \kappa_0 < \min a$\enskip regular
cardinals.  By [Sh:g, VIII \S 2] or
[Sh:506, \S2] $a$ has a generating sequence
$\l b_\tau \mid \tau \in pcf a\r$.  For
each $\xi < \theta$ we can find a successor
ordinal $\gamma_\xi \le \gamma^*$ so that
$\kap_*^{+\gamma_\xi} \in pcf (\{
\tau^k_\xi \mid k < \omega\}) \backslash \{
\lam_\ell\mid 1 \le \ell \le n^*\}$.  So,
for some successor ordinal $\gamma^{**} \le
\gam^*$ there is an unbounded in $\theta$ set
$Y$ consisting of $\xi$'s such that $\xi <
\theta$ and $\gamma_\xi = \gamma^{**}$.
Clearly, $\lam_\ell \in pcf a$ for $\ell =
1 \nek n^*$ and $\kap_*^{+\gam^{**}} \in
pcf a$.  Then w.l. of g. we can assume that
$b_{\kap_*^{+\gam^{**}}}$ is disjoint to
each $b_{\lam_\ell}$ for $\ell = 1\nek
n^*$. Set $A = \{ n < \omega\mid
b_{\kap_*^{+\gam^{**}}}\cap \sig_n$ is
unbounded in $\sig_n\}$.  
\medskip

\noindent
{\bf Claim 2.2} $A \in D$.
\medskip

\noindent
{\bf Proof.} If this does not hold, then
there is $\xi(*) < \theta$ such that for every
$n \in \omega \backslash A$
$b_{\kappa_*^{+\gam^{**}}}\cap
[\sig_{\xi(*)}, \sig_n) = \emptyset$.  W.l.
of g. $\xi(*) \in Y$.  Also, $n \in \omega
\backslash A$ implies that $\{
\sig^k_{\xi(*),n} \mid k < \omega\} \cap
b_{\kap_*^{+\gam^{**}}}=\emptyset$, since
for every $k < \omega$, $\sigma_{\xi(*)} <
\sig_{\xi(*),n}^k < \sig_n$.
\medskip

Hence $\{ \sig^k_{\xi(*),n}\mid k < \omega,
n \in \omega\backslash A\}$ is disjoint to
$b_{\kap_*^{\gam^{**}}}$.  Now, each
$\tau^k_{\xi(*)} \in pcf
(\{\sig^{k'}_{\xi(*),n} \mid k'<\omega, n
\in \omega\backslash A\})$.  Here we use
the assumption that $A \not\in D$ and so
$\omega \backslash A \in D$.

But $\kap_*^{+\gam^{**}} \in pcf (\{
\tau^k_{\xi(*)}\mid k < \omega\})$. Hence 
$\k \in pcf (\{ \sig^k_{\xi(*),n} \mid k <
\omega, n \in \omega\backslash A\})
\subseteq pcf (a \backslash b_\k)$,  which
is impossible by the choice of
generators.\\
$\square$ of the claim.

Let $n \in A$.  Then $b_\k \cap \sig_n$ is
unbounded in $\sig_n$.  Hence $pcf (b_\k
\cap \sig_n) \backslash \sig_n \ne
\emptyset$.  But $pp(\sig_n) \le
\sig_n^{+n^*}$, hence for some $\ell(n) \in
\{1\nek n^*\}$ we have $\sig_n^{+\ell(n)}
\in pcf (b_\k \cap \sig_n) \subseteq pcf
(b_\k)$.  Then for some $\ell(*) \in \{
1\nek n^*\}$ the set $A^* = \{ n \in A\mid
\ell(n) = \ell (*)\}$ belongs to $D$.  So,
$\lam_{\ell(*)} \in pcf (\{
\sig_n^{+\ell(*)} \mid n \in A^*\})
\subseteq pcf (b_\k)$. But $b_\k \cap
b_{\lam_{\ell(*)}} = \emptyset$.
Contradiction.\\
$\square$
\medskip

Using (3) of Theorem 2 we can give another
proof of Theorem 1.

\subsection*{2.3  Second proof of Theorem 1}
W.l. of g.  $cf\kap_* = (2^{\aleph_0})^+$.
Let $\theta = (2^{\aleph_0})^+$ and $\kap_0
= \theta^+$. Assume also w.l. of g. that $D$ is
a nonprincipal ultrafilter on $\omega$.
For every $f:\omega \to Reg \cap
\kap_1\backslash \kap_0$ we define
$g_f:\omega \to \rho^*<\theta$ as follows:
$$g_f (n) = \rho \quad {\rm iff}\quad f(n)
\in I_\rho\ .$$
Then, $f_1 \ge_D f_2$ will imply 
$g_{f_1} \ge_D g_{f_2}$ since the sequence
$\l \min I_\rho \mid \rho < \rho^*\r$ is
strictly increasing.  Consider $\l
f_{\lam_\alpha} \mid \alp < \theta\r$ of
(d) of Theorem 1.  This is a strictly
increasing sequence modulo $D$.  Now, the
total number of $g_f$'s is
$(\rho^*)^{\aleph_0} \le
(2^{\aleph_0})^{\aleph_0} = 2^{\aleph_0}$.
Hence there are $g^*: \omega \to \rho^*$ and
$\alp^* < \theta$ such that for every
$\alp, \theta > \alpha \ge \alp^*$, every $f:
\omega \to Reg \cap \kap_1\backslash
\kap_0$ such that $f_{\lam_\alp} \le_D f
<_D f_{\lam_{\alp+1}}$
$$f(n) \in I_{g^*(n)}, \ \hbox{for almost
each}\ n < \omega\ {\rm mod}\ D\ .$$
Apply $(*)_3$ to $\l g^*(n)\mid n <
\omega\r$ with $\gamma^* = 2$.  Then for
some $\ell^* \in \{1,2\}$ the following
holds:\\
if $a \in [\{ tcf (\prod_{n < \omega}
\sig_n /D)\mid \sig_n \in
I_{g^*(n)}$ for $n < \omega\} \cap
[\kap_1,\kap^*)]^{\aleph_0}$ then
$\kap_*^{+\ell^*} \not\in pcf a$. Let
$\alp$, $\theta > \alpha \ge \alp^*$.  Pick
$\kap_\alp$, $\lam_\alp < \kap_\alp <
\lam_{\alp + 1}$, $cf \kap_\alp =
\aleph_0$ and $pp\kap_\alp \ge \kap_*^{++}$
(by (c) of Theorem 1 we can assume w.l. of
g. that it exists). Then, by [Sh-g], there
are $\tau_{\alp,n} \in Reg \cap \kap_\alp
\backslash \lam_\alp^{++}$ $(n < \omega)$
and a filter $D_\alpha$ on $\omega$
containing all cofinite sets such that
$\kap_*^{+\ell^*} = tcf (\prod_{n < \omega}
\tau_{\alp,n} /D_\alp)$. Consider $\l
f_{\tau_{\alp,n}} (m) \mid m < \omega\r$
for every $n < \omega$. It is a sequence
of regular cardinals such that
$\tau_{\alp,n} = tcf (\prod_{m < \omega}
f_{\tau_{\alp,n}} (m) / D)$ and
$f_{\lam_\alp} <_D f_{\tau_{\alp,n}} <_D
f_{\lam_{\alp+1}}$.  Then for almost every
$m<\omega$ (mod $D$) $f_{\tau_{\alp,n}} (m)
\in I_{g^*(m)}$.  Hence $\tau_{\alp,n} \in
\{ tcf (\prod_{m < \omega} \sig_m / D)\mid
\sigma_m \in I_{g^*(m)}, m < \omega\}$ for
every $n < \omega$. Take $a =
\{\tau_{\alp,n}\mid n < \omega\}$. Then
$k_*^{+\ell^*} \not\in pcf a$, but
$\kap_*^{+\ell^*} = tcf (\prod_{n<\omega}
\tau_{\alp, n}/D_n)$. Contradiction.
\\
$\square$

\medskip
The following is parallel to 1.6.
\medskip

\noindent
{\bf Theorem 2.4}
{\it Suppose that
\begin{itemize}
\item[{\rm (a)}] $\kap_0 < \kap_1 < \kap_*$
\item[{\rm (b)}] $\theta_1,\theta_2 < \kap_0$ are
such 
that $cf \theta_1 > \aleph_0$, $\theta_2 =
\theta_1^{+3}$ or $\theta_2$ is regular
$\ge \theta_1^{+3}$ and for every $\alpha<
\theta_2$ $cf([\alp]^{<\theta_1},
\supseteq) < \theta_2$
\item[{\rm (c)}] $cf \kap_1 = \aleph_0$ and $pp\kap_1
\ge \kap_*^{+\theta_2}$
\item[{\rm (d)}] $\theta_3$ is regular cardinal between
$\theta_2$ and $\kap_0$
\item[{\rm (e)}] $\theta_4$ is cardinal between
$\theta_3$ and $\kap_0$ of cofinality $\ge
\theta_3$
\item[{\rm (f)}] $\theta_5 \in [\theta_4,\kap_0)$ is a
cardinal such that $cf ([\theta_5]^{\le
\aleph_0}, \subseteq) = \theta_5$
\item[{\rm (g)}] $D$ is an $\aleph_1$-complete filter on
$\theta_4 + 1$ 
\\
(Notice that we allow  $D$ to be
principal.  For example, generated by
$\{\theta_4\}$).
\item[{\rm (h)}] if $\l \mu_\alp \mid \alp \le
\theta_4\r$ is  a strictly increasing
continuous sequence of singular cardinals
between $\kap_0$ and $\kap_1$, then
$$\{\alp \le \theta_4 \mid \alp \
{\rm limit,}\ cf \mu_\alp \ge \theta_4
\ {\rm and}\ pp(\mu_\alp) <
\mu_\alp^{+\theta_1}\}\in D\ .$$
(Thus, if $\{\theta_4\}\in D$ then the
condition means $pp(\mu) < \mu^{+\theta_1}$
for every limit cardinal $\mu \in
(\kap_0,\kap_1)$ of cofinality $\theta_4$.)
\\
\underline{Then}
\item[{\rm (1)}] For every sequence $\overline\sig^* = \l
\sig^*_n \mid n < \omega\r$ of limit
cardinals of cofinality $\ge \theta_4$
between $\kap_0^+$ and $\kap_1$  there are
$\beta < \theta_2$ and a sequence
$\overline\sig^{**} = \l \sig_n^{**} \mid n <
\omega \r$, $\kap_0^+ \le \sig^{**}_n <
\sig^*_n$ $(n < \omega)$ such that\\[.5em]
$(\tilde*)_1$ if $a \in
[R_{\overline\sig^*,\overline\sig^{**}}]^{\aleph_0}$ then
$pcf(a) \cap [\kap_*^{+\beta},
\kap_*^{+\theta_2}) = \emptyset$, where
$R_{\overline\sig^*,\overline\sig^{**}} = \{\tau \in
(\kap_0^+, \kap_1)\mid$ there is a sequence
$\l\sig_n\mid n < \omega\r$, with $\sig_n
\in Reg \cap [\sig_n^{**}, \sig^*_n)$ such
that $\tau \in pcf \{ \sig_n\mid n <
\omega\}\}$.
\item[{\rm (2)}] There are $\alpha^* \le \theta_5$ and a
sequence $\l R_\alp \mid \alp < \alp^*\r$
with $\bigcup_{\alp < \alp^*} R_\alp = Reg
\cap \kap_* \backslash \kap_1$ so that\\[.5em]
$(\tilde *)_2$ for every $\alp < \alp^*$
there is $\beta < \theta_2$ such that for
every $a \in [R_\alp]^{\aleph_0}$ we have
$pcf (a) \cap [\kap_*^{+\beta},
\kap_*^{+\theta_2}) = \emptyset$.
\item[{\rm (3)}] There are $\rho^* < \theta^+_5$ and a
partition $\l I_\rho \mid \rho < \rho^*\r$
of $Reg \cap \kap_1\backslash \kap_0$ into
closed open intervals (i.e. of the form
$[x,y))$ with $\l \min I_\rho \mid \rho <
\rho^*\r$ strictly increasing such that
\\[.5em]
$(\tilde *)_3$ for every sequence of
ordinals $\l \rho_n\mid n < \omega \r$
below $\rho^*$ there is $\beta < \theta_2$
such that for every $a \in [\{tcf(\prod_{n
< \omega} \sig_n /\widetilde D)\mid \sig_n
\in I_{\rho_n}$ for $n < \omega$,
$\widetilde D$ is a nonprincipal
ultrafilter on $\omega$ with
$\lim_{n < \omega}{}_D (\min I_{\rho_n}) =
\kappa_1 \} ]^{\aleph_0}$
$$pcf (a) \cap
[\kappa_*^{+\beta} , \kap_*^{+\theta_2}) =
\emptyset\ .$$
\end{itemize}
}

\medskip
\noindent
{\bf Proof of (2) and (3) from (1)}

Let $\chi$ be a large enough cardinal.
Pick $M \prec (H(\chi), \epsilon)$ so
that $|M| = \theta_5,
\kappa_0,\kap_1,\theta_5 \in M$, $M \cap
\theta_5^+ \in \theta_5^+$ and $(\forall X
\in [M]^{\aleph_0}) (\exists Y \in M) (X
\subseteq Y \wedge |Y| = \aleph_0)$.

This is possible since by (f) $cf
([\theta_5]^{\le \aleph_0}, \subseteq) =
\theta_5$.  Define the set $\Phi$ now to be
$\{\overline\sig^* \in M\mid \overline\sig^*
= \l \sig_n^* \mid n < \omega\r$ is a
sequence of limit cardinals between
$\kap_0$ and $\kap_1$ with $cf \sig_n^* \ge
\theta_4$ $(n < \omega)\}$.

For each $\overline \sig^* \in \Phi$ we
choose
$\overline\sig^{**}=\overline\sigma^{**}
[\overline \sig^*]$ in $M$ satisfying
$(\tilde *)_1$. Define $\l R_\alp \mid \alp
< \alp^*\r$ to be an enumeration of the set
$\{R_{\overline\sig^*,\overline\sig^{**}[\overline\sig^*]} \mid
\overline\sig^* \in \Phi\} \cup \{ pcf
(\{\sig_n\mid n < \omega\})\mid \l \sig_n
\mid n < \omega\r \in \Phi$ and for every
$n < \omega\ \ cf \sig_n = \sig_n\}$.
 
Now we proceed as in Theorem 2. 
 
\medskip
\noindent
{\bf Proof of (1).}
Assume toward contradiction that for some
$\overline\sig^*$ there is no $\overline\sig^{**}$
satisfying (1).  We choose by induction on
$\xi < \theta_4$ cardinals
$\sig_{\xi,n}, \tau_\xi^{i,k},
\sig_{\xi,n}^{i,k}$ $(k,n<\omega,i <
\theta_2)$ such that
\begin{itemize}
\item[$(\alp)$] $\kap_0^+ \le
\sig_{\xi,n}<\sig_n^*$
\item[$(\beta)$] $\xi < \xi'$ implies that
$\sig^i_{\xi,n} < \sig_{\xi',n}$
\item[$(\gam)$] $\tau_\xi^{i,k} \in Reg \cap
\kap_*\backslash \kap_1$
\item[$(\delta)$] $pcf (\{
\tau_\xi^{i,k}\mid k < \omega\})\cap
[\kap_*^{+1}, \kap_*^{+\theta_2}) =
\emptyset$
\item[$(\epsilon$)] $\tau_\xi^{i,k} \in pcf
(\{ \sig_{\xi,n}^{i,k} \mid n < \omega\})$
\item[$(\xi)$] $\sig_{\xi,n} <
\sig^{i,k}_{\xi,n} = cf \sig^{i,k}_{\xi,n}
< \sig^*_n$
\item[$(\eta)$] $\l \sig_{\xi,n}\mid \xi <
\theta_4\r$ is an increasing continuous
sequence of singular cardinals.
\end{itemize}

The verification that such a construction is
possible is as in the proof of (1) of
Theorem 2.

Let $\sig_n = \sig_{n,\theta_4} =
\bigcup_{\xi < \theta_4} \sig_{\xi,n}$ for
each $n < \omega$.  Applying the condition
(h) of the statement of the theorem to $\l
\sig_{\xi,n} \mid \xi \le \theta_4\r$ we
find for every $n < \omega$ a set $Y_n \in
D$ such that $\xi \in Y_n$ implies that $pp
(\sig_{\xi,n}) < \sig_{\xi,n}^{+\theta_1}$.
By $\aleph_1$-completeness of $D$, the set
$Y = \bigcap_{n < \omega} Y_n \in D$.
Choose some $\delta^* \in Y$.  Let $pp
(\sig_{\delta^*,n}) =
(\sig_{\delta^*,n})^{+\beta_n}$ for some
$\beta_n < \theta_1$\enskip ($n < \omega$).

Consider sets $a_n = \{ \sig_{\xi,n}^{i,k}
\mid \xi < \delta^*, i < \theta_2, k <
\omega\}$ and $a = (\bigcup_{n < \omega}
a_n) \cup a^*$, where $a^* = \{
(\sigma_{\delta^*,n} )^{+\beta} \mid n <
\omega, \beta \le \beta_n$ is a successor
ordinal $\}$.  Then $a$ is a set of regular
cardinals of cardinality $\le \theta_4 +
\theta_2 < \kappa_0 < \min a$.  Let $\l
b_\tau \mid \tau \in pcf a\r$ be a
generating sequence.  As each $\beta_n <
\theta_1$ and $cf \theta_1 > \aleph_0$,
$|a^*| < \theta_1$.  By [Sh:g, IX] or
[Sh:g, E12, 4.18(b)] $c = pcf (a^*) \cap
[\kappa_*, \kap_*^{+\theta_2})$ is bounded
in $\kappa_*^{+\theta_2}$, since $\theta_2
\ge \theta_1 ^{+3} \ge |a^* |^{+4}$. Also
$pcf (c) = c$.  For each $\xi < \delta^*$
for some $i(\xi)$ we have $pcf (\{
\tau_\xi^{i(\xi),k} |k < \omega\}) \cap
[\kap_*, \kap_*^{+\theta_2})$ is not
bounded by $\sup c$.
So, choose $\kap_*^{+\rho (\xi)} \in pcf
(\{\tau_\xi^{i(\xi), k} \mid k < \omega\})
\cap [\kap_*, \kap_*^{+\theta_2})
\backslash \sup c$. Clearly, $\rho (\xi) <
\theta_2$ is a successor ordinal. As,
$\theta_2 < \theta_3 = cf \theta_3$, and
$\delta^*\in Y$ implies either $(cf\delta^*
= \theta_3)$ or $(\delta^* = \theta_4$ and
then also $cf \delta^*_1 \ge
\theta_3)$, necessary, for some $\rho^* <
\theta_2$ the set $Z = \{ \xi < \delta^*
\mid \rho(\xi) = \rho^*\}$ is unbounded in
$\delta^*$.  Let $J_n = J_{a_n}^{bd}$. So
$J_n$ is an ideal on $a_n$ and, clearly,
for every $c_n \in J_n$ $(n < \omega)$ we
have $\kap_*^{+\rho^*} \in pcf (\bigcup_{n
< \omega} (a_n \backslash c_n))$.

By $pcf$ theory (see [Sh:g, VIII, 1.5] or
[Sh:g, E12]) there are finite sets $e_n
\subseteq \cap \{ pcf (a_n \backslash
c_n)\mid c_n \in J_n\}$ $(n < \omega)$ such
that $\kap_*^{+\rho^*} \in pcf (\bigcup_{n <
\omega} e_n)$.  But $\cap \{ pcf (a_n
\backslash c_n) \mid c_n \in J_n\}
\subseteq \{\sig_{\delta^*,n}^{+\beta} \mid
\beta < \beta_n$ is a successor ordinal$\}$
for every $n < \omega$.  So $\bigcup_{n <
\omega} e_n \subseteq \cup \{
\sig_{\delta^*,n}^{+\beta} \mid \beta <
\beta_n$ is a successor ordinal and $n <
\omega\} = a^*$.  Hence, $\kappa_*^{+\rho^*}
\in pcf (a^*)$. But then $\kappa_*^{+\rho^*}
\in pcf (a^*) \cap [\kappa_*,
\kap_*^{+\theta_2}) = c$, which is
impossible by the choice of $\rho^*$.
Contradiction.\\
$\square$
\medskip

Let us conclude with a question which is
most natural taking into account the
results above.
\medskip

\noindent{\bf Question.} Is the following
situation possible:
\begin{itemize}
\item[(a)] $\kap_1 < \kap_*$, $cf \kap_1 =
\aleph_0$, $cf \kap_*=\aleph_1$
\item[(b)] for every singular $\mu <
\kap_1$, $pp\mu = \mu^+$ (or if one likes
only for $\mu$'s of countable cofinality)
\item[(c)] $\kap_* = \sup\{ \mu\mid \mu <
\kap_*, cf \mu = \aleph_0$ and $pp\mu =
\kap_*^+\}$
\item[(d)] the same as (d) of Theorem 1 or
even  add  $(*)$ of the proof of Theorem
1.
\end{itemize}

\end{document}